\documentclass[11pt]{amsart}
\usepackage{amssymb,amscd}
\usepackage[english]{babel}
\usepackage[latin1]{inputenc}
\usepackage[OT1]{fontenc}
\usepackage{times}
\usepackage{amsmath}
\usepackage{amsthm}
\usepackage{amsfonts}
\usepackage{amssymb}
\usepackage[all]{xy}
\newtheorem{theorem}{Theorem}[section]

\newtheorem{proposition}[theorem]{Proposition}

\theoremstyle{remark}
\newtheorem{remark}[theorem]{Remark}
\newtheorem{definition}[theorem]{Definition}

\newtheorem{example}[theorem]{Example}

\newcommand{\Lieg}{\mathfrak{g}}
\newcommand{\Lieh}{\mathfrak{h}}
\newcommand{\Liez}{\mathfrak{z}}

\newcommand{\np}{\mathfrak{n}_{+}}
\newcommand{\nm}{\mathfrak{n}_{-}}
\newcommand{\Wg}{W(\mathfrak{g})}
\newcommand{\Wh}{W(\mathfrak{h})}
\newcommand{\qWg}{\mathcal{W}(\mathfrak{g})}
\newcommand{\qWh}{\mathcal{W}(\mathfrak{h})}
\newcommand{\Sg}{S\mathfrak{g}}
\newcommand{\Sh}{S\mathfrak{h}}

\newcommand{\Ug}{U\mathfrak{g}}
\newcommand{\Sgg}{(S\mathfrak{g})^{\mathfrak{g}}}

\newcommand{\Zg}{Z(\mathfrak{g})}

\newcommand{\extg}{\bigwedge \Lieg}
\newcommand{\exth}{\bigwedge \Lieh}

\newcommand{\Cg}{\text{Cl}(\mathfrak{g},B_\Lieg)}

\newcommand{\Clg}{\text{Cl}\mathfrak{g}}
\newcommand{\Clh}{\text{Cl}\mathfrak{h}}

\newcommand{\C}{\mathbf{C}}
\newcommand{\corp}{\mathbf{k}}

\newcommand{\pg}{\widetilde{p}}

\newcommand{\qp}{\widehat{p}}

\newcommand{\HC}{\Phi_{HC}}
\newcommand{\Chev}{\Phi_{Ch}}

\begin{document}

\title{Transgression and Clifford algebras}

\author{Rudolf Philippe Rohr}
\address{University of Geneva, Departement of Mathematics,
2-4 rue du Lièvre, c.p. 64, 1211 Genève 4, Switzerland}
\email{Rudolf.Rohr@math.unige.ch}

\date{\today}

\begin{abstract}
Let $W$ be a differential (not necessarily commutative) algebra which carries a free action of a polynomial algebra $SP$ with homogeneous generators $p_1, \dots, p_r$. We show that for $W$ acyclic, the cohomology of the quotient $H(W/<p_1, \dots, p_r>)$ is isomorphic to  a Clifford algebra $\text{Cl}(P,B)$, where the (possibly degenerate) bilinear form $B$ depends on $W$. This observation is an analogue of an old result of Borel \cite{B} in a non-commutative context. As an application, we study the case of $W$ given by the quantized Weil algebra $\qWg = \Ug \otimes \Clg$ for $\Lieg$ a reductive Lie algebra. The resulting cohomology of the canonical Weil differential gives a Clifford algebra, but the bilinear form vanishes on the space of primitive invariants of the semi-simple part. As an application, we consider the deformed Weil differential (following Freed, Hopkins and Teleman \cite{FHT}).
\end{abstract}

\maketitle
\tableofcontents

\section{Introduction}

Let $G$ be a topological group, $EG$ the total space of the  universal bundle and $BG=EG/G$ the classifying space,

$$G \hookrightarrow EG \rightarrow BG.$$

In 1953 Borel \cite{B} showed how to compute the cohomology of $BG$ using the fact that the cohomology of $EG$ is trivial and 
the knowledge of the cohomology of $G$. More precisely, if $G$ has a homotopy type of a finite CW-complex, and with appropriate assumptions 
on the field $\corp$, we have 

$$H(G;\corp) \cong \bigwedge <x_1, ... , x_r>_{\corp},$$ 
as algebras. Then, using Borel's result, we have

$$H(BG;\corp) \cong \corp[y_1, ... ,y_r],$$
where the $y_i$'s are the images of the $x_i$'s under the transgression map (for more details see \cite{MC}, Theorem 3.27, Proposition 6.37 and Theorem 6.38).\\

In this article we are interested in the algebraic counterpart of Borel's argument. Let

$$SP \hookrightarrow W \longrightarrow F := W / <S^+P>,$$
where 

\begin{itemize}
	\item $(W,d)$ is a $\mathbf{Z}_2$-graded differential algebra;
	\item $SP$ is a symmetric subalgebra of $W$ generated by even central coboundaries $\{p_1,..., p_r \}$ such that $SP = \corp[p_1, ... ,p_r] \subset Z^{even}(W) \subset W$;
	\item $F := W / <S^+P>$ is the quotient of $W$ by the ideal generated by $p_1, \dots, p_r$.
\end{itemize}

Our main theorem (see §\ref{section_2}) states that if $W$ is a free $SP$-module and the cohomology of $W$ is trivial, we have 
an algebra isomorphism

$$H(F) \cong \text{Cl}(P,B).$$ 

Explicitly, this isomorphism is given by $p_i \rightarrow [C_{p_i}]_d \in H(F)$, where $C_{p_i} \in W$ is a cochain of transgression, $dC_{p_i} = p_i$. Even though the $C_{p_i}$'s are defined up to coboundary, the cohomology class $[C_{p_i}]_d$ is uniquely determined by $p_i$. The bilinear form $B$ is given by 

$$B(p_i,p_j) = \frac{1}{2}[C_{p_i} \cdot C_{p_j} + C_{p_j} \cdot C_{p_i}] \in H(W) \cong \corp.$$

In §\ref{section_3} we apply this theorem to the classical and quantized Weil algebras. Let $\Lieg$ be a complex reductive Lie algebra, and $B_{g}$ be a non degenerate invariant bilinear form on $\Lieg$. The classical Weil algebra is given by $\Wg = \Sg^* \otimes \bigwedge \Lieg^*$ and the quantized Weil algebra by $\qWg = \Ug \otimes \Cg$. Using results of Chevalley and Kostant, we have the following: $\Wg$ is a free $\Sgg$-module and $\Sgg = \mathbf{C}[\pg_1, ... , \pg_r]$ is generated by $r = \text{rank}(\Lieg)$ homogeneous polynomials, $\qWg$ is a free $\Zg$-module and $\Zg = \mathbf{C}[\qp_1, ... ,\qp_r]$, where $\qp_i$ is the image of $\pg_i$ under the Duflo isomorphism. We obtain the following result (Theorem \ref{thm_weil_algebra}):

$$H(\Wg / <(\Sgg)^+>) \cong \bigwedge <\pg_1, ... , \pg_r>_{\corp},$$
$$H(\qWg / <(\Zg)^+>) \cong \text{Cl} (<\qp_{1}, \dots , \qp_l>_\corp, B_{\Lieg}) \otimes \bigwedge <\qp_{l+1}, \dots , \qp_r>_\corp,$$
where $\{ \qp_{1}, \dots , \qp_l \}$ is a basis of the center of $\Lieg$, and $\{ \qp_{l+1}, \dots , \qp_r \}$ are primitive invariants of the semi-simple part of $\Lieg$.\\

We remark that on the semi-simple part, the cohomology is super-commutative.\\

In §\ref{section_4}, we follow Freed, Hopkins and Teleman \cite{FHT} to introduce a deformed Weil differential
$d' = d + \iota(\xi)$, where $d$ is the standard Weil differential and $\iota(\xi)$ is the contraction by an element $\xi$ of a Cartan subalgebra $\Lieh$ of $\Lieg$. On the $\Lieh$-invariant parts of $\Wg$ and $\qWg$, $d'$ is again a differential, and the cohomology in the classical and quantized cases is trivial. We obtain the following results (Theorem \ref{thm_def_W} and \ref{thm_def_qW} ):

$$H(\Wg^{\Lieh} / <(\Sgg)^+)> \cong \bigwedge P,$$
$$H(\qWg^{\Lieh} / <(\Zg)^+)> \cong \text{Cl} (P,B),$$
where the bilinear form $B$ is given, in terms of a basis $\{e_a\}_a$ of $\Lieg$, by 

$$B(p_i,p_j) = \sum_{a,b} \frac{\partial p_i}{\partial e_a}(\xi) \frac{\partial p_j}{\partial e_b}(\xi) B_\Lieg(e_a,e_b).$$

\textbf{Acknowledgments} I would like to thank Anton Alekseev, my thesis advisor, to have suggested me this problem and for very helpful discussions. I also tank Thierry Vust for very helpful discussions. This work was supported in part by the Swiss National Science Foundation.

\section{Cohomology of $W / <S^+P>$}
\label{section_2}

Throughout this section, $(W,d)$ denotes a $\mathbf{Z}_2$-graded differential algebra over a field $\corp$ of characteristic zero. We assume that its cohomology is trivial, i.e. $H(W,d) \cong \corp$.\\

Consider the super-center $Z(W) \subset W$, and its even part $Z^{even}(W) \subset Z(W) \subset W$. We choose $r$ linearly independent elements $\{p_1, \dots , p_r\} \subset Z^{even}(W)$ such that

\begin{enumerate}
	\item the $p_i$'s are coboundaries, i.e. $\forall i$ there exist $C_{p_i} \in W$ such that $dC_{p_i} = {p_i}$;
	\item the subalgebra of $W$ generated by the $p_i$'s and the unit is a symmetric algebra $SP$ of the vector space $P = <p_1, \dots ,p_r >_\corp$, i.e. $SP = \corp[p_1, \dots ,p_r]$.
\end{enumerate}

Let $S^+P$ be the augmented ideal of $SP$, i.e. $S^+P=\ker(\pi)$, where $\pi : SP \rightarrow \corp$ is the augmentation map sending all generators $\{p_i\}_{i=1 \dots r}$ to zero. Note that the differential vanishes on $SP$ and $S^+P$.\\

Let $<S^+P>$ be the ideal in $W$ generated by $S^+P$. As the differential vanishes on $S^+P$, it descends to the quotient algebra $W / <S^+P>$, and we obtain a new $\mathbf{Z}_2$-graded differential algebra $(W / <S^+P>,d)$. The main result of this section is the computation of its cohomology.\\

\begin{remark}
	The elements $C_p$ are defined up to coboundaries, but the classes $[C_p]_d \in H(W/<S^+P>,d)$ are uniquely determined by $p$. We call $C_p$ a cochain of transgression. This name comes from the classical Weil algebra, see \cite{K2}, Section 6.5.
\end{remark} 

\begin{remark}
We can begin with a $\mathbf{Z}$-graded algebra, by considering the induced $\mathbf{Z}_2$-grading.\\
\end{remark}

\begin{example}
Consider a finite dimensional vector space $V$, and its Koszul algebra $K(V) = SV \otimes \bigwedge V$ with its differential $d_{K(V)}$ (given on generators by $d_{K(V)}(1 \otimes x) = x \otimes 1$ and $d_{K(V)}(x \otimes 1) = 0$, see Chapter 3 of \cite{GS}). It is well known that the Koszul complex is acyclic. We have $SV \subset Z^{even}(K(V))$, and $SV = \corp[v_1, \dots , v_n]$ where $\{v_1 , \dots , v_n\}$ is a basis of $V$. Here we consider $P=V$. Then the quotient algebra $K(V)/<S^+V> = \bigwedge V$ is isomorphic to the exterior algebra and $d_{K(V)/<S^+V>} = 0$. This implies that $H(K(V)/<S^+V>) = \bigwedge V$.\\  
\end{example}

\begin{remark}
As $SP$ is contained in the even part of the super-center of $W$, $W$ has the structure of a $\mathbf{Z}_2$-graded $SP$-module. The action of $SP$ on $W$ is given by the multiplication on $W$.
\end{remark}

We remark that in the previous example the cohomology is given by the exterior algebra of the vector space $V$. In general, we have the following theorem:

\begin{theorem}
\label{cohomology_theorem}
	Let $W$ and $SP$ be as above. If $W$ is a free $SP$-module, we have an algebra isomorphism  
	
	$$H(W / <S^+P>, d ) \cong \text{Cl}\left(P,B \right).$$\\
	Moreover, this isomorphism is given by $P \ni p \rightarrow [C_p]_d \in H(W / <S^+P>)$, where $C_p$ is a cochain of transgression for $p$. The symmetric bilinear form $B$ is given by
	\begin{equation}
	\label{def_bilinear_form}
		\begin{array}{rcl}
			B:P \times P & \rightarrow & \corp, \\
			(p,q) & \mapsto & B(p,q) := \frac{1}{2}  \Big[ C_p \cdot C_q+C_q \cdot C_p \Big]_d.
		\end{array}
	\end{equation}
\end{theorem}

\begin{remark}
	The bilinear form is well defined;
	\begin{enumerate}
		\item for $p,q \in P$, we have $d[C_p,C_q] = 0$. This implies that $\Big[ [C_p,C_q] \Big]_d \in H(W,d) \cong \corp$;
		\item this definition does not depend on the choice of $C_p$ and $C_q$. Indeed, let $C'_p = C_p + da$ and $C'_q = C_q + db$, then $[C'_p,C'_q] = [C_p,C_q] + [C_p,db] + [da,C_q] + [da,db] = [C_p,C_q] + d[C_p,b] + d[a,C_q] + d([a,db])$.\\
	\end{enumerate}
\end{remark}

\begin{remark}
	The bilinear form vanishes if the algebra is $\mathbf{Z}_+$-graded \footnote{by $\mathbf{Z}_+$-graded we mean $\mathbf{Z}$-graded with all components of degree less than 0 vanishing}. Indeed, the elements $C_p$ are of degree at least $1$, hence $[C_p,C_q]$ is of degree at least $2$.\\
\end{remark}

\begin{example}
\label{example_class}
As in the previous example, we consider the Koszul algebra $K(V)$. But now, let $\mathcal{W}$ be a finite reflection group acting on $V$. Following Chevalley (see \cite{C}), the invariant part $(SV)^\mathcal{W}$ is generated by $n = dim(V)$ algebraically independent homogeneous polynomials, i.e. $(SV)^\mathcal{W} = \corp[p_1, ... ,p_n]$. Moreover $SV$ is a free $(SV)^{\mathcal{W}}$-module (see \cite{B1}, Chapter 5, Section 5.2). By applying Theorem \ref{cohomology_theorem} with $P = <p_1 , ... ,p_n >_\corp$, we get $H(K(V) / <((SV)^{\mathcal{W}})^+>  \cong \bigwedge P$ (the bilinear form vanishes because the Koszul algebra is $\mathbf{Z}_+$-graded), and a possible choice of cochains of transgression is given by
\begin{equation}
\label{equ_cochaine}
p \rightarrow \frac{1}{\text{deg}(p)+1} \frac{\partial p}{\partial e_a} \otimes e_a, \footnote{here we use the convention that we sum over a repeated index, i.e. in this case $\frac{\partial p}{\partial e_a} \otimes e_a = \sum_a \frac{\partial p}{\partial e_a} \otimes e_a$ }
\end{equation}
where $\{e_a\}$ is a basis of $V$.
\end{example}

\begin{example}
\label{example_quan}
As in the previous example, let $V$ be a vector space, and $\{p_1,..., p_n\}$ the generators of $(SV)^{\mathcal{W}}$. Here we consider a {\em deformed } version of the Koszul algebra, $\mathcal{K}(V,B_V) = SV \otimes \text{Cl}(V,B_V)$, with $B_V$ some bilinear form on $V$ (not necessarily non degenerate). The differential is given on generators by $d(1 \otimes x) = x \otimes 1$, whence $d(x \otimes 1) = 0$. The cohomology of $\mathcal{K}(V,B_V)$ is trivial. As in the previous example $\mathcal{K}(V,B)$ is a free $(SV)^{\mathcal{W}}$-module, and then with $P = <p_1 ,...,p_n >_\corp$, we get $H(\mathcal{K}(V,B_V) / <S^+P>) \cong \text{Cl}(P,B_V)$. Moreover we have the same choice for cochains of transgression as in the previous example. This allows us to calculate the bilinear form $B$:
\begin{equation}
\label{equ_bi}
B(p_i,p_j) = \Big[\frac{\partial p_i}{\partial e_a} \frac{\partial p_j}{\partial e_b} B_V(e_a,e_b) \Big]_d.
\end{equation}
Note that the bilinear form $B$ vanishes, if each polynomial is of degree at least $2$.
\end{example}

All the statements of these two examples will be proved in the next section.\\

The proof of Theorem \ref{cohomology_theorem} will proceed by induction and will use the universal-coefficient theorem for the cohomology of a $\mathbf{Z}_2$-graded modules (see Appendix \ref{appendix_1} for the proof). 

\begin{theorem}[Universal-coefficient theorem for the cohomology of $\mathbf{Z}_2$-graded modules]
\label{kunneth_formula}
Let $W$ be a $\mathbf{Z}_2$-graded free $R$-differential module, where $R$ is a principal commutative unitary ring, and let $M$ be a $R$-module. Then we have the following exact sequences:  
$$
\xymatrix{
0 \ar[r] & H^i(W) \otimes_R M \ar[r] & H^i(W \otimes_R M) \ar[r] & Tor_{R}(H^{i+1}(W),M) \ar[r] & 0
}
$$
where $i=0,1$, and these exact sequences split.\\

Moreover we have 

$$H^i(W \otimes_R M) = H^i(W) \otimes_R M \oplus (h^i)^*(Tor_{R}(H^{i+1}(W),M)),$$
where $(h^i)^*$ is the $R$-module homomorphism induced by the right inverse of the restriction of the differential to $W^i$.\\
\end{theorem}

To apply this theorem, we need the following $SP$-module.

\begin{definition}
Let $M = \mathbf{k}$ and define a $SP$-module structure by 

$$p \cdot v = \pi(p) \alpha,$$
where $\pi: SP \rightarrow \corp$ is the augmentation map.
\end{definition}

With this definition we have 

$$(W / <S^+P> ,d) \cong (W \otimes_{SP} M,d \otimes 1).$$\\

To prove Theorem \ref{cohomology_theorem}, we proceed by induction on the dimension of $P$, i.e. the number $r$ of generators $\{p_i\}_{i=1,...,r}$ of $SP$. For this we introduce the following sequence $\{(W_i,d_i) \}_{i=0,..., r}$ of differential algebras :

$$(W_0,d_0)=(W,d)\quad \text{and} \quad (W_i,d_i) = (W_{i-1} \otimes_{SP_{i}} M, d_{i-1} \otimes 1) \quad  r \geq i \geq 1,$$
where $P_i = <p_i>_{\corp}$. It is obvious that $W \otimes_{SP} M \cong W_r$, and if $W$ is a free $SP$-module, then $W_i$ is a free $SP_j$-module for all $j > i$.\\

We now take the first step, the calculation of the cohomology of $(W_1,d_1)$.\\

\begin{proposition}
\label{prop_one_element}
If $W_1$ is a free $SP_1$-module, we have $H(W_1,d_1) \cong  \text{Cl}\left(P_1,B \right)$. Moreover the isomorphism is given by $ p_1 \rightarrow [C_{p_1}]_{d_1} \in H(W_1,d_1)$, where $dC_{p_1} = p_1$.\\
\end{proposition}

\begin{proof}
\begin{enumerate}
\item
We apply the universal-coefficient theorem \ref{kunneth_formula} with $W=W_1$, $R = SP_1$ and $M$ as above. For $i=0$ we have the exact sequence
$$
\xymatrix{
0 \ar[r] & H^0(W_1) \otimes_{SP_1} M \ar[r] & H^0(W_1 \otimes_{SP_1} M) \ar[r] & 0,
}
$$
and for $i=1$  we have  
$$
\xymatrix{
0 \ar[r] & H^1(W_1 \otimes_{SP_1} M) \ar[r] & Tor_{SP_1} (H^0(W_1),M) \ar[r] & 0.
}
$$

These two exact sequences imply that $H(W_1) \cong \corp \oplus Tor_{SP_1} (H^0(W_1),M)$. We have to calculate the torsion product $Tor_{SP_1} (H^0(W_1),M)$. For this consider the following free resolution of $SP_1$-module over $H^0(W_1) = \corp$,
$$
\xymatrix{
0 \ar[r] & SP_1 \cdot p_1 \ar@{^{(}->}[r]^\iota & SP_1 \cdot 1 \ar[r] & \corp \ar[r] & 0, 
}
$$
where $\iota$ is the canonical injection. Now we make the tensor product with $M$. We then obtain the exact sequence of $SP_1$-modules

$$
\xymatrix{
0 \ar[r]& Tor_{SP_1} (\corp,M)  \ar[r] & SP_1 \cdot p_1 \otimes_{SP_1} M \ar[r]^{\iota \otimes 1} & SP_1 \cdot 1 \otimes_{SP_1} M \ar[r] & \corp \otimes_{SP_1} M \ar[r] & 0. 
}
$$
We have $Tor_{SP_1} (H^0(W_1),M) \cong \text{Ker}(\iota \otimes 1) = \mathbf{k} \cdot (p_1 \otimes 1)$. This implies that

$$H(W_1) \cong \corp \oplus \corp \cdot p_1,$$

as $\corp$-vector spaces.\\

\item
To give the isomorphism we need to calculate $(h^1 \otimes 1)(p_1 \otimes 1))$, where $h^1$ is a right inverse of $d_1$. Since $dC_{p_1} = p_1$, we can choose $(h^1 \otimes 1) (p_1 \otimes 1)) = (C_{p_1} \otimes 1)$, and then
 
$$(h^1)^*(Tor_{SP_1} (H^0(W_1),M)) = [C_{p_1} \otimes 1]_{d_1} \cong [C_{p_1}]_{d_1} \in H(W_1,d_1).$$\\
 
\item  
The last step is to compute the algebra structure. since $d_1$ is a derivation on $W_1$, the algebra structure descends to the cohomology, i.e. $[C_{p_1}]_{d_1} \cdot [C_{p_1} ]_{d_1} = [C_{p_1} \cdot C_{p_1}]_{d_1}$. Using the definition of the bilinear form (\ref{def_bilinear_form}), we have $[C_{p_1} \cdot C_{p_1} ]_{d_{1}} = B(p_1,p_1)$.\\

\end{enumerate}
\end{proof}

We now begin the proof of the main theorem.\\

\begin{proof}[Proof of Theorem \ref{cohomology_theorem}]
We proceed by induction on the dimension $r$ of $P$. For $r=1$ it is given by the above proposition. \\

\begin{enumerate}
	\item We use Theorem \ref{kunneth_formula} with $R = SP_{n+1}$ and $W = W_n$. We obtain the following exact split sequences ($i=0,1$):
$$\xymatrix{ 
0 \ar[r] & H^i(W_n) \otimes_{SP_{n+1}} M \ar[r] & H^i(W_{n+1}) \ar[r] & Tor_{SP_{n+1}}(H^{i+1}(W_n),M) \ar[r] & 0
}.$$
Then we have
$$H^i(W_{n+1}) = \text{Cl}^i( <[C_{p_1}],...,[C_{p_n}]>_\corp,B) \oplus (h^i)^* Tor_{SP_{n+1}}(H^{i+1}(W_n),M)$$
as $SP_{n+1}$-modules. The next step consists in calculating  the last term of this equality.\\

	\item Let $\{[C_{p_I}]\}$ be a $\mathbf{k}$-basis of $H^{i+1}(W_{n})$ ($I=(i_1, ... i_k)$ with $1 \leq i_1 < ... < i_k \leq n$ and $C_{p_{I}} = C_{p_{i_1}} \cdot ... \cdot C_{p_{i_k}}$). We obtain the following free exact sequence of $SP_{n+1}$-modules over $H^{i+1}(W_{n})$ :

$$
\xymatrix{
0 \ar[r] & \bigoplus_{I} SP_{n+1} \cdot p_{n+1}C_{p_{I}} \ar@{^{(}->}[r]^\iota & \bigoplus_{I} SP_{n+1} \cdot C_{p_{I}} \ar[r] & H^{i+1}(W_{n}) \ar[r] & 0,
}$$

and consequently $Tor_{SP_1} (H^{i+1}(W_1),M) = \bigoplus_{I} \mathbf{k} \cdot ( p_{n+1} C_{p_{I}} \otimes 1)$.\\

We have $d^i_{n+1}(C_{p_{n+1}} C_{p_{I}}) = p_{n+1} C_{p_{I}} $. This implies that $(h^i)^*Tor_{SP_{n+1}}(H^{i+1}(W_n),M) = \bigoplus_{I} \mathbf{k} \cdot ([C_{p_{n+1}}] [C_{p_{I}}])$, and then finally  

$$H^i(W_{n+1}) = \text{Cl}^i( <[C_{p_1}], ... ,[C_{p_{n}}]>_{\mathbf{k}},B) \oplus \bigoplus_{I} \mathbf{k} \cdot ([C_{p_{n+1}}]\cdot [C_{p_{I}}])$$ 

as $SP_{n+1}$-modules.\\

	\item The last step is the algebra structure. We have for all $k \leq n$ that $\Big[ [C_{p_{n+1}}, C_{p_k}] \Big ] = B(p_{n+1},p_k)$. This implies that $H(W_{n+1})$ injects in $\text{Cl}( <[C_{p_1}], ... ,[C_{p_{n+1}}]>_{\mathbf{k}},B)$. Since they have the same dimension, they are equal.\\
\end{enumerate}

\end{proof}

\begin{remark}
This proof of Theorem \ref{cohomology_theorem} still applies  if we assume only that $W$ is a free $SP_i$-module for all $i$.
\end{remark}


We now give an isomorphism theorem. Let $W_{I}$ and $W_{II}$ be two $\mathbf{Z}_2$-graded differentials. Denote by $SP_{I}$ and $SP_{II}$ the choice of the subalgebra of coboundaries elements in the even part of their center.

\begin{theorem}
\label{homo_theorem}
	Let $\Phi : W_I \rightarrow W_{II}$ be a graded differential vector space homomorphism, such that:
	\begin{enumerate}
		\item its restriction to $SP_I$ is an algebra isomorphism between $SP_{I}$ and $SP_{II}$, i.e. $\Phi : SP_I \stackrel{\cong}{\rightarrow} SP_{II}$;
		\item it commutes with the $SP$-module structure, i.e. we have $\Phi(p z) = \Phi(p) \Phi(z)$, $\forall p \in SP_I$ and $\forall z \in W_{I}$.
	\end{enumerate}
	
		 Then it induces a vector space isomorphism in cohomology,
	
	$$\bar{\Phi} : H(W_{I} / <S^{+}P_{I}>) \stackrel{\cong}{\rightarrow} H(W_{II} / <S^{+}P_{II}>).$$

	If in addition $\Phi$ is an algebra homomorphism, then this is an algebra isomorphism in cohomology.

\end{theorem}
\begin{proof}
		Firstly we remark that $\Phi$ induces a $SP$-module homomorphism 
		$$\bar{\Phi} : H(W_{I} / <S^{+}P_{I}>) \rightarrow H(W_{II} / <S^{+}P_{II}>)$$
		in cohomology. Secondly let $p_{II}$ be a generator of $SP_{II}$. Then there exists $p_I \in P_{I}$ such that $\Phi(p_I) = p_{II}$. Let $C_{p_I}$ be a cochain of transgression for $p_{I}$, then $\Phi(C_{p_I})$ is a cochain of transgression for $p_{II}$, i.e. $\bar{\Phi}$ is surjective and hence bijective.\\
		
 		For the last statement, we remark that if $\Phi$ is an algebra homomorphism, then $B(\Phi(p),\Phi(q)) = \Phi(B(p,q))$ for all $p,q \in SP_{I}$.\\
\end{proof}

\begin{remark}
	In order for the homomorphism $\Phi$ to induce an algebra isomorphism in cohomology, it suffices that the bilinear forms be isomorphic. 
\end{remark}

\section{Classical and quantized Weil algebras}
\label{section_3}

In this section, $\Lieg$ denotes a complex reductive Lie algebra. The classical Weil algebra $\Wg$ of $\Lieg$ is a $\mathbf{Z}_+$-graded $G$-differential algebra (see \cite{AM1}). We recall its definition and some elementary facts. The classical Weil algebra is defined by

$$\Wg = \Sg^* \otimes \bigwedge \Lieg^*.$$

The grading is given by degree $2$ on generators of $\Sg^*$ and degree $1$ on generators of $\bigwedge \Lieg^*$. Using an invariant non degenerate bilinear form we identify $\Lieg$ with its dual $\Lieg \cong \Lieg^*$.\\

The $G$-differential algebra structure is given by the following three derivations:\\

{\em the contraction} $\Lieg \ni x \rightarrow \iota(x) \in Der^{-1}(\Wg)$, where $\iota(x)(a \otimes b) = a \otimes \iota(x)b$ is the usual contraction on the exterior algebra, \\ 

{\em the adjoint action} $\Lieg \ni x \rightarrow L(x) \in Der^0(\Wg)$, where $L(x)(a \otimes b) = L_{\Sg}(x) a \otimes b + a \otimes L_{\bigwedge \Lieg }(x)b$ is given by the usual adjoint action on the symmetric and on the exterior algebra, i.e. the extension by derivation of the Lie bracket,\\

{\em the differential} $d \in Der^{1}(\Wg)$ is defined on generators by 

$$d(1 \otimes x) = x \otimes 1 + 1 \otimes \lambda(x),$$
where $\lambda: \Lieg^* \rightarrow \bigwedge^2 \Lieg^*$ is the dual of the Lie bracket.\\

These three derivations satisfy the following relations:

\begin{align*}
&[L(x),L(y)]=L([x,y]), \qquad [L(x),\iota(y)]=\iota([x,y]), \\
& \text{\em the Cartan formula} \quad [\iota(x),d]=L(x),
\end{align*}

and all other brackets vanish.\\

The cohomology of the Weil algebra is trivial, i.e. $H(W,d) \cong \C$. \\

The quantized Weil algebra $\qWg$ is an interesting deformation of $\Wg$. It is also a G-differential algebra, but is only $\mathbf{Z}_2$-graded. It is defined by 

$$\qWg = \Ug \otimes \Cg,$$
with $B_\Lieg$ some non degenerate invariant bilinear form on $\Lieg$ (see \cite{AM1} or \cite{AM2}). Let $\{e_i\}$ and $\{e^i\}$ be a pair of dual bases of $\Lieg$, and $f^{abc}$ be the structure constants, i.e. $[e_a,e_b] = f^{abc}e_c$.\\

In the quantized case the three derivations are given by commutators. Let 

\begin{align*}
g_i & = -\frac{1}{2} f^{iab} e_a e_b, \\ 
\mathcal{D} & = e^a \otimes e_a -\frac{1}{6}f^{abc} 1 \otimes e_a e_b e_c,
\end{align*}  
then 

$$\iota(e_i) = ad(1 \otimes e_i) \text{,} \quad L(e_i) = ad(e_i \otimes 1 + 1 \otimes g_i ) \quad \text{and } d = ad(\mathcal{D}).$$\\

There is an isomorphism of $G$-differential algebras between the classical and quantized Weil algebras, namely the quantization map $\mathcal{Q} : \Wg \rightarrow \qWg$. It is a $\mathbf{Z}_2$-graded vector space isomorphism which commutes with the contraction, the adjoint action and the differential (see Section 6 of \cite{AM1} or Section 4.3 of \cite{AM2}).\\

\begin{remark}
It is well known that there exists an algebra isomorphism between $\Sgg$ and the center $\Zg$ of $\Ug$, the Duflo isomorphism (see \cite{Du}). The quantization map restricts to the Duflo isomorphism, on $\Sgg \otimes 1$.
\end{remark}

\begin{remark}
The restriction of the quantization map to the exterior algebra $ 1 \otimes \bigwedge \Lieg $ gives the usual Chevalley symmetrization map.
\end{remark}

Let $\Lieh$ be a Cartan subalgebra of $\Lieg$; we can consider Weil algebras $\Wh$ and $\qWh$. In these two cases, the adjoint action vanishes and the differential is less complicated. We have as differential algebra that they are the Koszul algebras of $\Lieh$. More precisely we have $\Wh = K(\Lieh)$ and $\qWh = \mathcal{K}(\Lieh,B_\Lieh)$, where $B_\Lieh$ is the restriction to $\Lieh$ of $B_\Lieg$. Moreover for the classical Weil algebra there exists a differential algebra homomorphism. This is the Chevalley projection

$$\Chev : \Wg \rightarrow \Wh = K(\Lieh),$$\\
which is the restriction homomorphism. For the quantized case there is a differential space homomorphism, the Harish-Chandra projection

$$\HC : \qWg \rightarrow \qWh = \mathcal{K}(\Lieh,B_\Lieh).$$\\

These two homomorphisms and the Duflo isomorphism give the following commutative diagram of graded algebra isomorphisms:

\begin{equation}
	\xymatrix{
		\Zg  \ar[rd]_{\HC} & \Sgg \ar[l]_{\text{Duf}} \ar[d]^{\Chev} \\
		& S\Lieh^{\mathcal{W}} }
\end{equation}
where $\mathcal{W}$ is the Weyl group of $\Lieg$. Using the same generators of $(S \Lieh )^{\mathcal{W}}$ as in Example \ref{example_quan}, we obtain

\begin{align*}
\Sgg &= \C[\pg_1, ... \pg_r], \\
\Zg & = \C[\qp_1, ... \qp_r],
\end{align*}
with the following relations:  
  
  $$\HC(\qp_i) = \pg_i \text{,} \quad \text{Duf}(\pg_i) = (\qp_i) \text{,} \quad \Chev(\pg_i) = p_i.$$\\
  
See Appendix \ref{appendix_2} for more details about the Harish-Chandra and Chevalley projections.\\

We give a description of these invariant polynomials. Let $\Lieg = \Liez \oplus \Lieg'$ be the decomposition of $\Lieg$ into its center $\Liez$ and its semi-simple part $\Lieg'$. On the center we have $S(\Liez)^\Liez = \C[\Liez]$. Then we can choose $\{\pg_1, \dots, \pg_l \}$ to be a basis of $\Liez$, where $l = \dim(\Liez)$. This gives 

$$\Sgg = \C[\pg_1, \dots, \pg_l, \pg_{l+1}, \dots \pg_r],$$
where $\{\pg_{l+1}, \dots \pg_r \}$ are the generators on the semi-simple part. Moreover we can choose the polynomial $\pg_{l+1}$  to be the Casimir polynomial, and $\text{deg}(p_{i}) = m_{i}+1$ where $1 = m_{l+1} \leq m_{l+2} \leq ... \leq m_{r}$. These integers are the exponents of the Lie algebra $\Lieg'$. \\

In the three cases elements of $\Sh^W$, $\Sgg$ and $\Zg$ are coboundaries. So a natural choice for the vector space $P$ for the Weil algebras is\\

{\em for Weil algebras}, $P_{\Wg} = <\pg_1, ... , \pg_r>_\C$,\\

{\em and for quantized Weil algebras}, $P_{\qWg} = <\qp_1, ... \qp_r>_\C$.\\

In \cite{K1}, Kostant proved that $\Sg$ is a free module over $\Sgg$ and that $\Ug$ is a free module over its center $\Zg$. This implies that $\Wg$ is a free module over $\Sgg$ and that $\qWg$ is a free module over $\Zg$. \\

With this choice it is obvious that the hypotheses of Theorem \ref{cohomology_theorem} are satisfied, and consequently the main results of this section are,

\begin{theorem}
\label{thm_weil_algebra}

We have, as algebras, that

\begin{enumerate}
	\item $H(\Wg / <S^+P_{\Wg}> ) \cong \bigwedge P_\Wg,$
	\item $H(\qWg / <S^+P_{\qWg}> ) \cong \text{Cl} (<\qp_{1}, \dots , \qp_l>_\C, B_{\Lieg}) \otimes \bigwedge <\qp_{l+1}, \dots , \qp_r>_\C.$
\end{enumerate}

\end{theorem}

\begin{remark}
	A choice of cochain of transgression for the elements of the semi-simple part of $P_{\Wg}$ is given in \cite{K2} Theorem 62. The quantizations of the cochains of transgression of $\Wg$ are cochains of transgression for $\qWg$,  i.e.  $\mathcal{Q}( C_{\pg} )= C_{\mathcal{Q}(\pg)}$.
\end{remark}

\begin{remark}
The algebra structure of the classical case is obvious. Indeed the classical Weil algebra is $\mathbf{Z}_+$-graded.
\end{remark}

\begin{remark}
Using Theorem \ref{homo_theorem}, we conclude that the Chevalley projection $\Chev \Wg \rightarrow \Wh = K(\Lieh)$ induces an isomorphism in cohomology, i.e. 
$$\Chev : H(\Wg/  <S^+P_{\Wg}>) \stackrel{\cong}{\rightarrow} H(K(\Lieh) / <S^+P>).$$ 
\end{remark}

To conclude the proof of the theorem we have to establish the algebra structure in the quantized case. For this we will prove that the Harish-Chandra projection 

$$\HC : \qWg \rightarrow \qWh = \mathcal{K}(\Lieh,B_{\Lieh})$$
induces an algebra isomorphism in cohomology. Then with the statements of Example \ref{example_quan}, with $\mathcal{W}$ the Weil group of $\Lieg$ and $V = \Lieh$, this finishes the proof. But first we will prove all statements of Example \ref{example_quan}. For this we will construct an explicit homotopy operator. The Koszul differential is given on generators by

$$d(1 \otimes x) = x \otimes 1 \qquad \text{and then} \qquad d(x \otimes 1) = 0.$$

Define a derivation $s$ on generators by 

$$s(1 \otimes x) = 0 \text{ and then } s(x \otimes 1) = x.$$ 

We have $[d,s] = Id$ on generators . We want to show that the inclusion of scalars $i : \C \rightarrow \mathcal{K}(V,B_V)$ and the augmentation map $\pi : \mathcal{K}(V,B_V) \rightarrow \C$ (defined by sending generators to zero) are homotopy inverse. Since $[d,s]$ is a derivation, $[d,s]+ i \cdot \pi$ is invertible on $\mathcal{K}(V,B_V)$, and we can define 

\begin{equation}
\label{homotopy_operator}
	h = s \cdot ([d,s]+ i \cdot \pi)^{-1}.
\end{equation}

This operator is the desired homotopy operator, i.e. we have $[d,h] = id - i \cdot \pi$.\\

Using the above homotopy operator on $p \in P$ we obtain $(d \cdot h )p = p$. This implies that a choice for cochains of transgression is given by

$$p \rightarrow h(p) = \frac{1}{\text{deg}(p)+1}\frac{\partial p}{\partial e_i} \otimes e_i \in SV \otimes \text{Cl}(V,B).$$

This establishes (\ref{equ_cochaine}). It is obvious that the bilinear form is given by (\ref{equ_bi}).\\

With $V=\Lieg$ the bilinear form vanishes on the semi-simple part. Indeed all polynomials of $\{p_{l+1}, \dots ,p_{r}\}$ are of degree at least two. On the center we have $B(p_i,p_j) = B_{\Lieg}(p_i,p_j)$ ($1 \leq i,j \leq l$).\\

The last step in the proof of Theorem \ref{thm_weil_algebra} is to prove that the Harish-Chandra projection induces an algebra isomorphism in cohomology.\\

\begin{proposition}
\label{prop_HC_Ch}
	The Harish-Chandra projection $\HC^ : \qWg \rightarrow \qWh = \mathcal{K}(\Lieh,V)$ satisfies the hypotheses of Theorem \ref{homo_theorem}. Moreover it induces an algebra isomorphism in cohomology.\\
\end{proposition}
\begin{proof}
	It satisfies hypothesis (a) because its restriction to $\Zg$ and $\Sh^W$ is an algebra isomorphism. Using (\ref{equ_hc_sp}) we conclude that it satisfies hypothesis (b). For the bilinear form, we have that cochains of transgression are in the $\Lieg$-invariant part of $\qWg$. But on the $\Lieh$-invariant part the Harish-Chandra projection is an algebra homomorphism. Then we can conclude that $B(p,q) = B(\HC(p),\HC(q))$ for all $p,q \in \qWg$. This implies the algebra isomorphism in cohomology.\\
\end{proof}

\section{Deformation of the Weil differential}
\label{section_4}

In \cite{FHT}, Freed, Hopkins and Teleman introduce a deformation of the Weil differential in the quantized case. For this, fix an element $\xi \in \Lieh$.  The deformation of $\mathcal{D}$ is then given by 

$$\mathcal{D'} = \mathcal{D} -  1 \otimes \xi.$$  

This provides a new even derivation $d'$ given by $d' = d - \iota{(\xi)}$. Its square is given by $d'^{2} = -[d, \iota{(\xi)}] = -L(\xi)$ (using Cartan's formula). Hence it is a differential on the kernel of $L({\xi})$.\\

The main result of this section is to give results analogous to Theorem \ref{thm_weil_algebra} and Example \ref{example_quan} for the deformed differential.\\
 
\subsection{The case of the Koszul algebra, $\mathcal{K}(V,B)$}

Let $\{e_a\}$ be a basis of $V$. Decompose $\xi$ in this basis, say $\xi = \xi^{a}e_{a}$. The deformed differential is given on generators by

$$d' (1 \otimes e_a) = - \xi^a (1 \otimes 1) + e_a \otimes 1 \quad \text{and} \quad d'(e_a \otimes 1) = 0.$$

Using the same derivation $s$ as in the non deformed case, we get

$$[d',s] (1 \otimes e_a) = (1 \otimes e_a) \quad  [d',s] ((e_a - \xi^a) \otimes 1) = (e_a - \xi^a) \otimes 1.$$

Then in the new variables $1 \otimes e_a$ and $\xi^a - e_a \otimes 1$ we arrive at the usual Koszul differential. We have, for the same homotopy operator $h$ as in (\ref{homotopy_operator}),

$$[d',h] = I - i \cdot \pi.$$

But the augmentation map $\pi$ is now defined as sending $(e_a - \xi^a) \otimes 1$ and $1 \otimes e_a$ to zero.\\

With this homotopy operator we conclude that the cohomology is trivial, i.e. $H(\mathcal{K}(V,B_V),d') \cong \corp$.\\

Let $\{p_1, \dots p_r\}$ be the generators of $S(V)^\mathcal{W}$ in Example \ref{example_quan}. With the deformed differential they are no longer coboundaries, since $(d' \cdot h  )p_i = p_i - p_i(\xi)$. But $p_i -p_i(\xi)$ are coboundaries.\\

Then a good choice for the vector space $P$ is  $P = <p_1-p_1(\xi), ... ,p_r-p_r(\xi)>_\corp$. The new ring $SP$ is isomorphic to the old one, thus $\mathcal{K}(V,B_V)$ is now $SP$-free. But note that the new ideal $S^+P$ is not isomorphic the old one. In the new case when we quotient we send $p_i - p_i(\xi)$ to zero instead of sending $p_i$ to zero.\\

We will now calculate the bilinear form (\ref{def_bilinear_form}). But first, with the homotopy operators, we calculate cochains of transgression . The invariant polynomials $p_i$ are given in the variables $e_a$. But to use the homotopy operator we make a change of variables by setting $u_a = e_a - \xi_a$. Using Taylor series we have
	
	$$p = p(\xi) + u_a \frac{\partial p}{\partial e_a} (\xi) +  \frac{1}{2}u_a u_b \frac{\partial^2 p}{\partial e_a \partial e_b} (\xi) + \dots  ,$$
whence for cochains of transgression,

	$$\C_p = e_a  \otimes \frac{\partial p}{\partial e_a}(\xi) + e_a \otimes  u_b \frac{\partial^2 p}{\partial e_a \partial e_b} (\xi) + \dots   .$$
		
This gives for the bilinear form

\begin{equation}
\label{equ_bi_2}
	B(p - p(\xi),q -q(\xi)) = \frac{\partial p}{\partial e_a}(\xi) \frac{\partial q}{\partial e_b}(\xi) B_V(e_a,e_b).
\end{equation}
	
Then using Theorem \ref{cohomology_theorem} we obtain the following result.\\
	
\begin{proposition}
\label{prop_koszul}
	The cohomology of $(\mathcal{K}(V,B_V),d')$ is given by
		$$H((\mathcal{K}(V,B_V),d') / <S^+P>) \cong \text{Cl}(P,B).$$
\end{proposition}

\subsection{The classical case}

In the classical case, the deformation of the Weil differential is again a differential if we restrict to the $\Lieh$-invariant part. For this reason we will restrict to $\Wg^{\Lieh}$. Using the same notations as in the previous section we have on a basis of $\Lieg$  

$$d'(1 \otimes e_a)  = e_a \otimes  - \xi^a - 1 \otimes \lambda(e_a).$$

This is more complicated than for the Koszul algebra. But there exists a change of variables which brings us back to a Koszul algebra. For more details see Chapters 3.1 and 3.2 of \cite{GS}.\\

In the new variables $z_a := e_a \otimes 1 - 1 \otimes \lambda(e_a)$ we have

$$d' z_a = - L(\xi)( 1 \otimes e_a) \quad \text{and} \quad d'(1 \otimes e_a )= z_a - \xi^a.$$ 

We will now give a homotopy operator to prove that the cohomology of $\Wg^{\Lieh}$ is trivial.\\

Define a derivation $s$ on generators by $s (1 \otimes e_a)= z_a$. We obtain  

$$[d',s](z_a-\xi^a) = z_a-\xi^a \quad \text{and} \quad [d',s](1 \otimes e_a) = 1 \otimes e_a.$$

Now as before, $h = s \cdot ([d',s] + i \cdot \pi)^{-1}$ is a homotopy operator, i.e. $[d',h] = I - i \cdot \pi$. But here the augmentation map $\pi$ sends $1 \otimes e_a$ and $z_a - \xi^a$ to zero. This implies that $H(\Wg^{\Lieh},d') \cong \C$.\\

Note that this augmentation map can also be defined by sending $e_a \otimes 1 - \xi^{a}$ and $1 \otimes e_a$ to zero.\\

Now as for the Koszul algebra, the $\pg_i$ are not coboundaries, but the $\pg_i - \pg_i(\xi)$ are. Indeed, we have $[d',h](\pg_i-\pg_i(\xi))= (d' \cdot h)(\pg_i-\pg_i(\xi)) = \pg_i-\pg_i(\xi)$. Then we choose $P_{\Wg} = < \pg_1-\pg_1(\xi),...,\pg_r-\pg_r(\xi)>_\corp$. It is obvious that $\Wg^{\Lieh}$ is $SP_{\Wg}$-free. Indeed we have $\Sg = SP_{\Wg} \otimes A$, which implies that $\Wg^{\Lieh} = SP_{\Wg} \otimes (A \otimes \bigwedge \Lieg)^{\Lieh}$.\\

Using Theorem \ref{cohomology_theorem} we obtain

\begin{theorem}
\label{thm_def_W}
	The cohomology of $(\Wg^{\Lieh} / <S_{\Wg}^+P> , d')$ is given by
	
	$$H(\Wg^{\Lieh} / <S_{\Wg}^+P>) \cong \bigwedge P_{\Wg}.$$
\end{theorem}

The bilinear form vanishes. Indeed, this can be seen either by using the $\mathbf{Z}_+$-grading for which the degree of $1 \otimes e_a$ is one and the degree of $e_a \otimes 1 - \xi^{a}$ is two, or the Chevalley projection which induces an algebra isomorphism in cohomology, i.e. $\bar{\Phi}_{Ch} : H(\Wg^{\Lieh} / <S_{\Wg}^+P>) \stackrel{\cong} {\rightarrow} H(K(\Lieh)) \cong \bigwedge P$.

\subsection{The quantized case}

In the quantized case the deformation of the Weil differential is again a derivation on the $\Lieh$-invariant part. The quantization map induces a $\mathbf{Z}_2$-graded vector space isomorphism in cohomology:

$$\mathcal{Q}: H(\Wg, d') \stackrel{\cong}{\rightarrow} H(\qWg, ad(\mathcal{D}')).$$

Using the results of the preceding subsection we have

\begin{align*}
&H(\qWg^{\Lieh},ad(\mathcal{D}')) \cong \C,\\
&\qp_i - \pg_i(\xi) \text{ are coboundaries}.
\end{align*}

Consequently we define the $SP_{\qWg}$-module structure by $P_{\qWg} = <\qp_1-\pg_1(\xi), ... \qp_r-\pg_r(\xi)>_\corp$. It is obvious that $\qWg^{\Lieh}$ is $SP_{\qWg}$-free. We remark that on the $\Lieh$-invariant part the Harish-Chandra projection is an algebra homomorphism. Hence it induces an algebra isomorphism in cohomology. Using Proposition \ref{prop_koszul} and Theorem \ref{cohomology_theorem} we have\\

\begin{theorem}
\label{thm_def_qW}
	We have, as algebras,
	
	$$H(\qWg^{\Lieh} / <S^+P_{\qWg}>, ad(\mathcal{D}')) \cong \text{Cl}(P_{\qWg},B),$$
	where $B$ is the same bilinear form as for the Koszul algebra, cf. (\ref{equ_bi_2}), with $p_i = \HC(\qp_i)$.
\end{theorem}

\begin{appendix}
\section{Universal-coefficient theorem for cohomology of $\mathbf{Z}_2$-graded modules}
\label{appendix_1}

In this appendix we prove Theorem \ref{kunneth_formula}. It is  a $\mathbf{Z}_2$-graded version of the universal theorem for homology of $\mathbf{Z}$-graded modules which can be found in \cite{S} Chapter 5, Section 2 \footnote{in \cite{S} the torsion product is denoted by $A \ast_R B = Tor_R(A,B)$}. We call it the {\em  Universal-coefficient theorem for cohomology of $\mathbf{Z}_2$-graded modules}, because in the $\mathbf{Z}_2$-graded cases there is no distinction between cohomology and homology.\\

Throughout this section $R$ denotes a commutative unitary ring, $W$ a differential $\mathbf{Z}_2$-graded $R$-module and $d$ its differential. We suppose that the grading of elements in the ring $R$ is even. Then we have $W = W^0 \oplus W^1$, where $W^0$ is the $R$-module of the even elements of $W$ and $W^1$ is the $R$-module of its odd elements, and the following diagram for the differential:\\

$$
\xymatrix{
W^1 \ar@/^/[r]^{d} &  W^0 \ar@/^/[l]^{d}
}.
$$

The proof of this theorem is in three parts. In the first and second parts we prove two exact sequences from which we deduce the exact sequence of the theorem. And in the last part we give a right inverse to $H^i(W \otimes_R M) \rightarrow Tor_R(H^i(W),M)$ and then prove that the exact sequence splits. Before we beginning the proof we will fix some notation and give some elementary facts. \\

Consider ($i=0,1$) the coboundaries spaces

$$B = Im(d) \quad \text{and} \quad W^i \supset B^i := B \cap W^i = \{d a | a \in W^{i-1} \},$$
and

$$Z = Ker(d) \quad \text{and} \quad W^i \supset Z^i := Z \cap W^i = \{a \in W^i| d a =0 \}.$$

Since the differential module $W$ is $\mathbf{Z}_2$-graded its cohomology inherits the grading. We have

$$H(W) = H^0(W) \oplus H^1(W) \text{, where } H^i(W) = Z^i / B^i.$$

Considering $B$ and $Z$ as differential $R$-modules with vanishing differential we have the exact sequences of differential $R$-modules ($i=0,1$)

\begin{equation}
\label{exact_seq_base}
\xymatrix{
0 \ar[r] & Z^i \ar@{^{(}->}[r]^{\alpha^i} & W^i \ar[r]^{d^i} & B^{i+1} \ar[r] & 0
},
\end{equation}
where $\alpha^{i}$ is the canonical injection and $d^i$ the restriction of $d$ to $W^{i}$. 


Let $M$ be a $R$-module. We can regard $W \otimes_R M$ as a differential $\mathbf{Z}_2$-graded $R$-module with differential $d \otimes 1$. Now we will consider the tensor product of the above exact sequence with $M$. If in the exact sequence (\ref{exact_seq_base}) the $R$-module $B^i$ is free then the differential $R$-modules sequence 
$$
\xymatrix{
0 \ar[r] & Z^i \otimes_R M \ar@{^{(}->}[r]^{\alpha^i \otimes 1} & W^i \otimes_R M \ar[r]^{{d^i} \otimes 1} & B^{i+1} \otimes_R M \ar[r] & 0 
}
$$
is exact (see Lemma 3.3 in Ch. XVI of \cite{L}). From this we obtain the exact {\em "ring"} of $R$-modules

\begin{equation}
\label{exact_ring}
\xymatrix{
& Z^{0} \otimes_R M \ar[r]^{(\alpha^0 \otimes 1)^*} & H^0(W \otimes_R M ) \ar[dr]^{(d^0 \otimes 1)^*} \\
B^0 \otimes_R M \ar@{^{(}->}[ur]^{q^{0} \otimes 1} & & & B^{1} \otimes_R M \ar@{^{(}->}[dl]^{q^1 \otimes 1} \\
& H^1(W \otimes_R M) \ar[ul]^{(d^1 \otimes 1)^*} & Z^1 \otimes_R M \ar[l]^{(\alpha^1 \otimes 1)^*}  
}
\end{equation}

It is a $\mathbf{Z}_2$-graded version of the usually long exact sequence which we meet in the $\mathbf{Z}$-graded case (see Section 1, Chapter 5 of \cite{S}).   \\

We conclude these preliminaries by giving the torsion product between $H(W)$ and $M$. For this, we consider the exact sequence of $R$-modules

$$
\xymatrix{
0 \ar[r] & B^i  \ar[r]^{q^i} & Z^i \ar@{->>}[r] & H^i(W) \ar[r] & 0
}.
$$

If $W$ is a free differential $R$-module, i.e. $W^0$ and $W^1$ are free $R$-modules, and $R$ is a principal ring, then the previous exact sequence is a free presentation of $H^i(W)$. Indeed $Z^i$ and $B^i$ are free $R$-modules (Corollary 3, Ch. VII.14 of \cite{B2}). Then by the characteristic property of the torsion product we have the exact sequence of $R$-modules

\begin{equation}
\label{second_exact_sequence}
\xymatrix{
0 \ar[r]& Tor_R(H^i(W),M) \ar[r] & B^i \otimes_A M  \ar[r]^{q^i \otimes 1} & Z^i \otimes_R M \ar@{->>}[r] & H^i(W) \otimes_R M \ar[r] & 0
}.
\end{equation}

Moreover we have $Tor_R(H^i(W),M) = Ker(q^i \otimes 1)$.\\

\begin{proof}[Proof of Theorem \ref{kunneth_formula}]
\begin{enumerate}
	\item From (\ref{exact_ring}) we get the following exact sequence \footnote{$\text{coker}(q^i \otimes 1) = (Z^i \otimes_R M) / Im(q^i \otimes 1)$}:
	
$$
\xymatrix{
0 \ar[r] & \text{coker}(q^i \otimes 1) \ar[r] & H^i(W \otimes_R M) \ar[r] & \text{ker}(q^{i+1} \otimes 1) \ar[r] & 0
}.
$$ 

	\item From the exact sequence (\ref{second_exact_sequence}), we have $\text{coker}(q^i \otimes 1) = H^i(W) \otimes_R M$ and $\text{ker}(q^i \otimes 1) = Tor_R(H^i(W),M)$. Thus we obtain 
$$
\xymatrix{
0 \ar[r] & H^i(W) \otimes_R M \ar[r] & H^i(W \otimes_R M) \ar[r] & Tor_R(H^{i+1}(W),M) \ar[r] & 0
},
$$
where $H^i(W) \otimes_R M \rightarrow H^i(W \otimes_R M)$ is the map induced by the bilinear map $ H^i(W) \times M \rightarrow H^i(W \otimes_R M)$, which assigns to $([w],m)$ the class $[w \otimes m]$, where $w$ is a cocycle of $W^{i}$ and $m \in M$.
	
	\item Thirdly we will give a right inverse to $H^i(W \otimes_R M) \rightarrow Tor_R(H^i(W),M)$. Since $B^{i+1}$ is a free $R$-module and $Im(d^i) = B^{i+1}$, there exists a $R$-module homomorphism $h^{i} : B^{i+1} \rightarrow W^{i}$ such that $d^i \circ h^{i} = Id$, i.e. which is a right inverse of $d^i$. Then the map 
	
$$h^i \otimes 1 : B^{i+1} \otimes_R M \rightarrow W^{i} \otimes_R M$$
sends $\text{ker}(q^{i+1} \otimes 1)$ into a cocycle of $W^{i} \otimes_R M$, and induces the desired map

\begin{equation}
\label{right_inverse}
	(h^i)^* : Tor_R(H^{i+1}(W),M) \rightarrow H^{i}(W \otimes_R M).
\end{equation}
Remark that the map $h_{i}$ is not unique but $h^*_i$ is unique.\\ 
	
\end{enumerate}
\end{proof}

\begin{remark}
The theorem still applies if we assume only that $R$ is a commutative unitary ring and $Z^i$ and $B^i$ are projective $R$-modules.
\end{remark}

\section{The projections of Harish-Chandra and of Chevalley}
\label{appendix_2}

In this appendix we recall the definition and some properties of the Harish-Chandra and Chevalley projections. First some notation: $\Lieg$ denotes a reductive Lie algebra, $\Lieh$ its Cartan subalgebra, $B_\Lieg$ some invariant non degenerate bilinear form. Let $\Lieg = \nm \oplus \Lieg \oplus \np$ be its triangular decomposition.\\

We begin with the Chevalley projection which we encountered in the classical case, and then we treat the Harish-Chandra projection which we met in the quantized case.

\subsection{Chevalley projections}

For the symmetric algebra of $\Lieg$ we have the decomposition $\Sg = \Sh \oplus (\nm \Sg + \Sh \np)$. This allow us to define an algebra homomorphism by the projection on the first term, i.e. $\Chev^{\Sg} :\Sg \rightarrow \Sh$. This is the original Chevalley projection. This projection can be viewed as the restriction homomorphism.\\

We recall Theorem 7.3.7 of \cite{D}, which says that it restricts to a graded algebra isomorphism $\Chev^{\Sg} : \Sgg \stackrel{\cong}{\rightarrow} \Sh^{W}$.\\

After the symmetric algebra we have the exterior algebra. We have the decomposition $\extg = \exth \oplus (\nm \extg + \extg \np)$. Then we define the Chevalley projection as the projection on the first term, i.e. $\Chev^{\extg} : \extg \rightarrow \exth$. It is an algebra homomorphism and can be viewed as the restriction homomorphism. But in this case the non trivial invariants are mapped to zeros (see \cite{BA}, Corollary 5.4.6).\\

For the Weil algebra we define $\tilde{\mathfrak{n}}_+ = 1 \otimes \np + \np \otimes 1$ and similarly for $\tilde{\mathfrak{n}}_-$. We have the decomposition $\Wg = \Wh \oplus (\tilde{\mathfrak{n}}_- \Wg + \tilde{\mathfrak{n}}_+ \Wg))$. This allows us to define the Chevalley projection $\Chev : \Wg \rightarrow \Wh$. It is clearly an algebra homomorphism and can be viewed as the restriction homomorphism. It is given by the tensor product of the Chevalley projections in the symmetric and exterior algebras, i.e. $\Chev = \Chev^{\Sg} \otimes \Chev^{\extg}$. Moreover it is a graded differential algebra homomorphism. \\

\subsection{Harish-Chandra projections}

Let $\Ug$ be the enveloping algebra of $\Lieg$. Following Chapter 7.4 of \cite{D}, we have the decomposition $\Ug^{\Lieh} = \Sh \oplus L$ where $L = \Ug \np \cap \Ug^\Lieh$ is a two sided ideal. This can by generalized by $\Ug = \Sh \oplus (\nm\Ug + \Ug\np)$. Then we have a projection $\kappa : \Ug \rightarrow \Sh$. Let $\gamma$ the automorphism of the $\Sh$ algebra which transforms the polynomial function $p$ into $\lambda \rightarrow p(\lambda - \rho)$, where $\rho$ is the half-sum of positive roots. The composition $\HC^{\Ug} = \gamma \circ \kappa$ is the Harish-Chandra projection for the enveloping algebra. Note that in general it is not an algebra homomorphism. But on the subspace $\Lieh$-invariant subspace it is.\\

We recall Theorem 7.4.5 of \cite{D}, which says that the Harish-Chandra projection restricts to an algebra isomorphism $\HC^{\Ug} : \Zg \stackrel{\cong}{\rightarrow} \Sh^{W}$. With the Chevalley projection and the Duflo isomorphism we obtain the following commutative diagram of graded algebra isomorphisms:

\begin{equation}
	\xymatrix{
		\Zg  \ar[rd]_{\HC^{\Ug}} & \Sgg \ar[l]_{\text{Duf}} \ar[d]^{\Chev^{\Sg}} \\
		& \Wh^{W} }
\end{equation}

Moreover we have $\forall p \in \Zg$ and $\forall q \in \Ug$ the identity

\begin{equation}
\label{equ_hc_sp}
	\HC^{Ug}(p \cdot q) = \HC^{Ug}(p) \cdot \HC^{Ug}(q).
\end{equation}

Indeed, let $p = \tilde{p} + \bar{p}$ be the decomposition of $p$, where $\tilde{p} = \HC^{Ug}(p)$ and $\bar{p} \in \nm\Ug + \Ug\np$. And let $q = n_-hn_+$ be the decomposition of $q$ in the PBW basis of $\Ug$, where $n_-$, $h$ and $n_+$ are products of respectively negative nilpotent elements, Cartan's subalgebra elements and positive nilpotent elements. In this notation we have $p \cdot q = n_-h\tilde{p}n_+ + n_-h\bar{p}n_+$. In order that the Harish-Chandra projection not vanish, we must have $n_- = n_+ =1$. In this case we obtain $\HC^{\Ug}(p \cdot q) = h \cdot \tilde{p} = \HC^{\Ug}(p)\cdot\HC^{\Ug}(q)$. In the other case, both sides vanish.\\

After the enveloping algebra we have the Clifford algebra. We have the decomposition $\Clg = \Clh \oplus (\nm \Clg + \Clg \np)$. This defines the Harish-Chandra projection $\HC^{\Clg} : \Clg \rightarrow \Clh$. The restriction to the $\Lieh$-invariant part gives an algebra homomorphism (see \cite{BA}, Chapter 5).\\

For the quantized Weil algebra, we have the decomposition  $\qWg = \qWh \oplus (\nm \qWg + \qWg \np)$. This defines the Harish-Chandra projection $\HC : \qWg \rightarrow \qWh$. Note that in general this is not an algebra homomorphism, but on the $\Lieh$ invariant part is it. As in the classical case, we have that $\HC = \HC^{\Ug} \otimes \HC^{\Clg}$ (see Section 7 of \cite{AM2}). Moreover it is a graded differential space homomorphism.\\

\end{appendix}


\bibliographystyle{amsplain}   %

\end{document}